\documentclass[11pt,notitlepage,twoside]{article}
\usepackage{latexsym} 
\usepackage{amsmath} 
\usepackage{amsfonts} 
\usepackage{amssymb} 
\renewcommand{\a }{\alpha } 
 
\renewcommand{\d}{\delta }

\newcommand{\D }{\Delta }
 
\newcommand{\e }{\varepsilon }

\renewcommand{\l }{\lambda } 
\renewcommand{\L }{\Lambda } 
 
\newcommand{\n }{\nabla }

\newcommand{\Sig }{\Sigma}

\renewcommand{\O }{\Omega } 
 
\newcommand{\ov}{\overline}

\newcommand{\be}{\begin{equation}} 
\newcommand{\ee}{\end{equation}} 
 
\newenvironment{pfn}[1]{\noindent{\bf Proof of {#1}\enspace}}{
\hfill$\Box$\medskip} 
\newcommand{\R}{\mathbb{R}} 
 
\newcommand{\N}{\mathbb{N}}

\newtheorem{thm}{Theorem}[section] 
\newtheorem{pro}[thm]{Proposition}
\newtheorem{lem}[thm]{Lemma}

\newtheorem{cor}[thm]{Corollary}
 
\numberwithin{equation}{section}

\author{Hichem CHTIOUI$^a$ \& Khalil EL MEHDI$^{b,c}$  \\ 
{\footnotesize
a : D{\'e}partement de Math{\'e}matiques, Facult{\'e} des Sciences de Sfax, Route
Soukra,}\\
{\footnotesize
 Sfax, Tunisia. E-mail : \texttt{Hichem.Chtioui@fss.rnu.tn} }\\
{\footnotesize
b :  Facult\'e des Sciences et Techniques, Universit\'e de Nouakchott, BP 5026, Nouakchott,}\\{\footnotesize
  Mauritania. E-mail : \texttt{khalil@univ-nkc.mr}}\\
{\footnotesize
 c : The Abdus Salam International Centre for Theoretical Physics, Mathematics Section,}\\
{\footnotesize
 Strada Costiera 11, 34014 Trieste, Italy. E-mail : \texttt{elmehdik@ictp.trieste.it}}
}

 \title { {\Large {\textbf { On a Paneitz Type Equation in Six Dimensional Domains}}}}

\begin{document}

\date{ }

\maketitle

{\footnotesize 
\noindent 
{\bf Abstract.}
In this paper, we consider the equation  $\D ^2 u=K u^5$, $u>0$ in $\Omega$, $u=\D u=0$ on $\partial\Omega$, where $K$ is a positive function and $\Omega$ is a bounded and smooth domain in $\R^6$. Using the theory of critical points at infinity, we give some topological conditions on $K$ to ensure some existence results.
}

\section{Introduction }
\mbox{}
 Let $\O$ be a smooth bounded domain of $\R^n$ with 
 $n\geq 5$ and consider the following nonlinear  problem
 under the Navier boundary condition
$$
(P)\qquad \left\{\begin{array}{ccccc}
\Delta ^2 u & = & K  u^p, & u > 0 &\mbox {in}\,  \Omega\\
u & = & \D u & =0 &\mbox {on}\,  \partial \Omega
\end{array}
\right.
$$
where $K$ is a $C^3$-positive function and where $p+1=\frac{2n}{n-4}$ is the critical exponent of the
embedding $H^2\cap H_0^1(\O)$ into $L^{p+1}(\O)$.\\
The interest in this problem grew up from its resemblance to some
nonlinear equation arising from a
geometric  context. Namely, the problem of prescribing the Paneitz curvature, which consists in finding suitable conditions on a given function $K$ defined on $M^n$ sucht that $K$  is the Paneitz curvature for a metric  $\tilde{g}$ conformally equivalent to $g$, where $(M^n, g)$ is a $n$-dimensional compact Riemannian manifold (for details on can see
\cite{BE1},\cite{BE2}, \cite{DHL},\cite{DMO}, \cite{F}
and the references therein).\\
Fourth order equations with critical growth on a domain of $\R^n$ has been studied in earlier works (see \cite{BEH}, \cite{BGP}, \cite{BH}, \cite{EFJ}, \cite{EO}, \cite{HV}, \cite{Lin}, \cite{NSJ},  \cite{PV},\cite{PS}, \cite{V1} and \cite{V2}).\\
The special nature of problem $(P)$ appears when we consider it from a variational viewpoint, indeed the Euler Lagrange functional associated to $(P)$ does not satisfy the Palais-Smaile condition, that is  theres exist noncompact sequences along which the functional is bounded and its gradient goes to zero. This fact is  due to the presence of the critical exponent. Moreover
it is easy to see that a necessary condition for solving the
problem $(P)$ is that $K$ has to be positive somewhere. In addition, there is at least another obstruction to solve problem $(P)$, based on
the Kazdan-Warner type condition, see \cite{DHL}. Hence it is not expectable to solve problem $(P)$ for all functions $K$, and it is natural to ask : under which conditions on $K$, $(P)$ has a solution?\\
Ben Ayed, El Mehdi and Hammami \cite{BEH} gave some sufficient topological conditions on $K$ to ensure the existence of solutions of $(P)$ for $n$ bigger than or equal to $6$. Ben Ayed and Hammami \cite{BH} provided an Euler-Hoph type criterion for $K$ to find solutions of $(P)$ for $n=6$. The methods of \cite{BEH} and \cite{BH} involve the study of critical points at infinity and their contribution to the topology of the level sets  of the associated Euler Lagrange functional.\\
Furthermore, as already known for problems related to the Scalar Curvature, there is a new phenomenon  in dimension $n\geq 7$, due the fact that the self interaction of the functions failing Palais-Smaile condition dominates the interaction of two of those functions. In the five dimensional case, the reverse happens. In the six dimensional case, we have a balance phenomenon, that is the self interaction and the interaction are of the same size (see \cite{DMO}, \cite{BE1}, \cite{BE2}, \cite{BEH}, \cite{BH}).\\ \
In this paper, we single out the six dimensional case to give more existence results. We are thus loking for  solutions of the following problem  
$$
(1)\qquad \left\{\begin{array}{ccccc}
\Delta ^2 u & = & K  u^5, & u > 0 &\mbox {in}\,  \Omega\\
u & = & \D u & =0 &\mbox {on}\,  \partial \Omega
\end{array}
\right.
$$
We precisely follow some of the ideas developed in Bahri \cite{B2}, Ben Ayed-Chtioui-Hammami \cite{BCH}, Chtioui \cite{C} and Chtioui-El Mehdi \cite{CE}. The main idea is to use the difference of topology of the critical points at infinity between the level sets of of the associated Euler Lagrange functional and the main issue is under our conditions on $K$, there remains some difference of topology which is not due to the critical points at infinity and therefore the existence of solutions of $(1)$.\\
In order to state our results, we need to fix some notations and  assumptions that we are using in our results.

We denote by $G$ the Green's function and by $H$ its regular part, that is for each $x\in \O$,
$$
 \left\{\begin{array}{ccccc}
G(x,y) & = & \mid x-y\mid^{-(n-4)}-H(x,y)  &  \mbox { in }\,  \O  \\
\D ^2 H(x,.)& = & 0  & \mbox { in }\, \O  \\
\D G(x,.) & = & G(x,.) =0 & \mbox{ on } \partial \O
\end{array}
\right.
$$
Throughout  this paper, we assume that  the following two assumptions hold\\
 ${\bf 1.1}$\hskip 0.3cm $K$ has only nondegenerate
critical points $y_0$, $y_1$,...,$y_m$ such that $y_0$ is the unique absolute maximum of $K$ on $\overline{\O}$ and such that
$$-\frac{1}{60}\frac{\D K(y_i)}{K(y_i)}+H(y_i,y_i)\ne 0 \quad \forall\,
i=0,1,...,m.$$\\ \\
 ${\bf 1.2}$\hskip 0.3cm All the critical points of
$K_1=K_{/\partial\O}$ are $z_1,...,z_{m'}$, and satisfy
$$
\frac{\partial K}{\partial \nu}(z_i) < 0, \qquad \mbox {for } i=1,...,m',
$$
where $\nu$ is the outward normal to $\O$.\\
We now introduce the following set
$$
\mathcal{F}^+ =\{ y_i\in \O/ \n K(y_i)=0 \,\, \mbox{ and }\,\,-\frac{1}{60}\frac{\D K(y_i)}{K(y_i)}+H(y_i,y_i) > 0 \}
$$
Thus our first result is the following :
\begin{thm}\label{t:11}
If $y_0 \notin \mathcal{F}^+$, then
  problem $(1)$ has a solution.
\end{thm}
To state our second result, we need to fix some notation.\\
For $s\in\N^*$ and for any $s$-tuple $\tau_s=(i_1,...,i_s) \in (\mathcal{F}^+)^s$ such that $i_p\neq i_q$ if $p\neq q$, we define a Matrix $M(\tau_s)=(M_{pq})_{1\leq p,q\leq s}$, by
$$
M_{pp}=\frac{-\D K(y_{i_p})}{60K(y_{i_p})^{\frac{3}{2}}}+\frac{
H(y_{i_p},y_{i_p})}{K(y_{i_p})^{\frac{1}{2}}}, \qquad M_{pq}
=-\frac{G(y_{i_p},y_{i_q})}{\left(K(y_{i_p})K(y_{i_q})\right)^{1/4}}
\quad \mbox{for } p\neq q,
$$
and we denote by $\rho (\tau _s)$ the least eigenvalue of $M(\tau _s)$.
It was first pointed out by Bahri \cite{B1} (see also \cite{BaC},
 \cite{BCCH} and \cite{R}), that when the self interaction and the
interaction between the functions failing Palais-Smale condition are the same size,
the function $\rho$ plays a fundamental role in the existence
of solutions to problems like $(P)$. Regarding problem $(P)$,
Ben Ayed and Hammami \cite{BH} observed that such
kind of phenomenon appears for $n=6$.\\
Now let $Z$ be a pseudogradient of $K$ of Morse-Smale 
type (that is the intersections of the stable and the unstable manifolds of the critical points of $K$ are transverse).\\
${\bf 1.3}$\hskip 0.3cm Assume that
 $W_s (y_i)\cap W_u(y_j) = \emptyset $ for any $y_i\in \mathcal{F}^+$ and for any $y_j\notin \mathcal{F}^+ $, where  $W_s (y_i)$ is the stable manifold of $y_i$ and $W_u(y_j)$ is the unstable manifold  of $y_j$ for $Z$.\\
Let $X=\bigcup_{y_i\in \mathcal{F}^+}\overline{W_s(y_i)}.$\\
${\bf (H_1)}$\hskip 0.3cm Assume that $X$ is not contractible.\\
We then have the following result
\begin{thm}\label{t:12}
Under the assumption   $(H_1)$, if  the following two conditions hold
\begin{align*}
(C_1)& \quad  \mbox{for any }\,s,\,\,M(\tau_s) \,\,\mbox {is nondegenerate}\\
(C_2)& \quad \rho (y_i,y_j) < 0 \quad \forall y_i,\, y_j \in \mathcal{F}^+ \,\, \mbox{such that}\,\, y_i\neq y_j,
\end{align*}

then problem $(1)$ has a solution.
\end{thm}
In Theorem \ref{t:11}, we have assumed that $y_0 \notin \mathcal{F}^+$. Next we want to give some existence result for problem $(1)$ when $y_0 \in \mathcal{F}^+.$ To this aim, we introduce some notation.\\
${\bf (H_2)}$\hskip 0.3cm Assume that $y_0 \in \mathcal{F}^+$\\
Let $y_{i_1} \in \mathcal{F}^+ \diagdown \{y_0\}$ such that\\
${\bf 1.4}$\hskip 0.3cm
$K(y_{i_1}) = \max \{K(y)/y\in \mathcal{F}^+ \diagdown \{y_0\}\}$\\
 and we denote by $k_1=6-i(y_{i_1})$, where  $i(y_{i_1})$ is the Morse index of $K$ at  $y_{i_1}$.\\
${\bf (H_3)}$\hskip 0.3cm Assume that $i(y_{i_1} \leq 5$.\\
We then have the following result
\begin{thm}\label{t:13}
 Under  assumptions $(H_2)$ and $(H_3)$, if the following three conditions hold
\begin{align*}
(A_0)& \quad  M(y_0,y_{i_1}) \,\,\mbox {is nondegenerate}\\
(A_1)& \quad \rho (y_0,y_{i_1}) < 0\\
(A_2)& \quad \frac{1}{K(y)^{1/2}} > \frac{1}{K(y_0)^{1/2}} + \frac{1}{K(y_{i_1})^{1/2}} \qquad \forall y \in \mathcal{F}^+ \diagdown \{y_0, y_{i_1}\},
\end{align*}
then problem $(1)$ has a solution of Morse index $k_1$ or $k_1+1$. 
\end{thm}
In contrast to Theorem \ref{t:13}, we have the following results based on a topological invariant for some Yamabe type problems introduced by Bahri \cite{B2}. To state these results, we need to fix  assumptions that we are using and some notation.\\
${\bf (H_4)}$\hskip 0.3cm Assume that 
$$
\biggl(-\frac{1}{60}\frac{\D K(y_0)}{K(y_0)}+H(y_0,y_0)\biggr)\biggl(-\frac{1}{60}\frac{\D K(y_{i_1})}{K(y_{i_1})}+H(y_{i_1},y_{i_1})\biggr) > G^2(y_0,y_{i_1}) $$
Let 
$$
X=\overline{W_s(y_{i_1})}
$$
Notice that under ${\bf 1.3}$ and ${\bf 1.4}$, we see that $X = W_s(y_{i_1})\cup W_s(y_0)$ and therefore $X$ is a manifold of dimension $k_1$ without boundary.\\
${\bf (H_5)}$\hskip 0.3cm Assume that $K(y_0) > 4K(y_{i_1})$ \\
Now we denote by $C_{y_0}(X)$ the following set 
$$
C_{y_0}(X)=\{\a \d _{y_0}+(1-\a )\d _x \, / \, \a \in [0,1],\, x\in X \},
$$
where $\d_x$ is the Dirac measure at $x$.\\
For $\l $ large enough, we introduce the map $f_\l : C_{y_0}(X)\to \Sigma ^+$,
defined by
$$
 (\a \d _{y_0}+(1-\a )\d _x)  \longrightarrow 
\frac{\a \d _{(y_0,\l )}
+(1-\a )\d _{(x,\l )}}{||\a \d _{(y_0,\l )}+(1-\a )\d _{(x,\l )}||},
$$
where
 $$
\Sig ^+=\{u\in H^2\cap H^1_0(\O )/ u\geq 0 \quad \mbox{and}\quad ||u||^2 =1\} \qquad \mbox{with} \quad ||u||^2:=\int _\O(\D u)^2
$$
and where $\d_{(x,\l )}$ is defined in the next section by \eqref{d}.\\
Then $ C_{y_0}(X)$ and $f_\l ( C_{y_0}(X))$ are manifolds in dimension $k_1+1$, 
that is, their singularities arise in dimension $k_1-1$ and lower, see \cite{B2}. We observe that 
 $ C_{y_0}(X)$ and $f_\l ( C_{y_0}(X))$ are contractible while $X$ is 
not contractible.\\
For $\l $ large enough, we also define the intersection number(modulo 2) of 
$f_\l (C_{y_0}(X))$ with $ W_s(y_0,y_{i_1})_\infty$
$$
\mu (y_0)=f_\l (C_{y_0}(X)). W_s(y_0,y_{i_1})_\infty ,
$$ 
where $W_s(y_0,y_{i_1})_\infty$ is the stable manifold of the critical point at infinity
$(y_0,y_{i_1})_\infty$ (see Corollary \ref{c:32} below) for a 
decreasing pseudogradient $V$ for the Euler Lagrange functional associated to $(1)$ which is transverse to 
$f_\l (C_{y_0}(X))$. Thus this number is well defined, see \cite{M}.\\
We then get the following result
\begin{thm}\label{t:14}
 Under assumptions  $(H_3)$, $(H_4)$ and $(H_5)$,
if $\mu (y_{i_1})=0$ then problem $(1)$ has a solution of Morse index $k_1$ or $k_1+1$.
\end{thm}
Now we give a statement  more general than Theorem \ref{t:14}. To this aim, let $k\geq 1$ and define $X$ as the following
$$
X=\overline{\cup_{y\in B_k} W_s(y)}, \quad \mbox{ with } B_k \,\,\mbox{ is any subset in }\,\,\{y \in \mathcal{F}^+ /
\, \mbox{i}(y)=6-k \},
$$
where $i(y)$ is the Morse index of $K$ at $y$.\\
${\bf (H_6)}$\hskip 0.3cm
We assume that $X$ is a stratified set without boundary (in the topological sense, that is, $X\in \mathcal{S}_k(S^4_+)$, the group of chains of dimension $k$ and $\partial X= \emptyset$).\\
${\bf (H_7)}$\hskip 0.3cm Assume that for any critical point $z$ of $K$ in $X\diagdown \{y_0\}$, we have 
$$
\biggl(-\frac{1}{60}\frac{\D K(y_0)}{K(y_0)}+H(y_0,y_0)\biggr)\biggl(-\frac{1}{60}\frac{\D K(z)}{K(z)}+H(z,z)\biggr) > G^2(y_0,z).
$$
For $y\in B_k$ we define, for $\l $ large enough, the intersection
number(modulo 2) 
$$
\mu (y)=f_\l (C_{y_0}(X)). W_s(y_0,y)_\infty
$$
By the above arguments, this number is well defined, see \cite{M}. \\
${\bf (H_8)}$\hskip 0.3cm Assume that $K(y_0) >4K(y)$ $\forall y\in \mathcal{F}^+ \diagdown \{y_0\}$.\\
Then we have the following theorem
\begin{thm}\label{t:15}
 Under  assumptions  $(H_6)$, $(H_7)$ and $(H_8)$,
if $\mu_(y)=0$ for any $y\in B_k$, then problem $(1)$ has a solution of Morse index $k$ or $k+1$.
\end{thm}

The plan of the rest of the  paper is the following. In section 2, we set up the variational structure and recall some preliminaries, while section 3 is devoted to  the proofs of our results. 

\section { Some Preliminaries }
\mbox{}
In this section we recall the functional setting and the variational problem associated to $(1)$. We will also recall some useful  results.\\
Problem $(1)$ has a variational structure, the related functional is
 $$
 J(u) =\biggl(\int_\Omega\ K \mid
u\mid^6\biggr)^{-\frac{1}{3}}
$$
 defined on
$$
\Sigma =\{ u\in H^2\cap H_0^1(\O) / \mid\mid u\mid\mid_{H^2\cap
H_0^1(\O)}^2:=\mid\mid u\mid\mid_2^2:= \int_\O\mid\D u\mid ^2 =1\}.
$$
The positive critical points of $J$ are solutions to (P), up a
multiplicative constant. Due to the non-compactness of the
embedding $H^2\cap H_0^1(\O)$ into $L^{p+1}(\O)$, the functional
$J$ does not satisfy the Palais-Smale condition. To characterize  the sequences failing the Palais-Smale condition, we need to fix some notation.\\
 For $a\in\Omega $  and $ \lambda > 0 $, we introduce the following function
\begin{eqnarray}\label{d}
\delta_{(a,\lambda )}(x)= c_0 \frac{\lambda} {1+\lambda
^2\mid x-a\mid^2},
\end{eqnarray}
$c_0$ is chosen so that $\delta_{(a,\lambda)}$ is the family of
solutions of the following problem (see \cite{Lin})
\begin{eqnarray}
\D ^2 u=\mid u \mid^{4}u, \quad u>0\quad \mbox{ in }
\R^6.
\end{eqnarray}
 For $f\in H^2(\O)$, we define the projection $P$ by:
\begin{eqnarray}
u=Pf\Longleftrightarrow\D ^2u=\D^2f \mbox{ in } \O,\quad u=\D u=0
\mbox{ on } \partial \O.
\end{eqnarray}
We now introduce the set of potential critical points at infinity.\\
For $\e>0$ and $p\in \N^*$, let $V(p,\e)$ be the subset of
$\Sig$ of the following functions: $u\in \Sig$ such that there is
$(a_1,...,a_p)\in \O^p$, $(\l_1,...,\l_p)\in (\e^{-1},+\infty)^p$
and $(\a_1,...,\a_p)\in (0,+\infty)^p$ such that
$$
\bigg |\bigg | u-\sum_{i=1}^p\a_i
P\d_{(a_i,\l_i)}\bigg |\bigg |_{2}<\e,\,\,
\l_id(a_i,\partial\O)>\e^{-1},\,\, \bigg |
\frac{\a_i^4 K(a_i)}{\a_j^{4}K(a_j)}-1\bigg | <\e,\,\,
\e_{ij}<\e \mbox{ for } i\ne j, 
$$
 where
\begin{eqnarray}
\e_{ij}=\biggl(\frac{\l_i}{\l_j}+\frac{\l_j}{\l_i}+\l_i\l_j\mid
a_i-a_j\mid^2\biggr)^{-1}.
\end{eqnarray}
The failure of the Palais-Smale condition can be described, following the ideas introduced in \cite{BrC},\cite{HR}, \cite{Lio}, \cite{S} as follows:
\begin{pro}\label{p:21}  
Assume that $J$ has no critical point in $\Sigma ^+$ and let $(u_k)\in 
\Sigma ^+$ be a sequence such that $J(u_k)$ is bounded and $\n J(u_k)\to 0$.
Then, there exist an integer $p\in \N^*$, a sequence $\e _k>0$ ($\e _k\to 0$)
and an extracted subsequence of $u_k$, again denoted $(u_k)$, such that 
$u_k\in V(p,\e _k )$.
\end{pro}
We now introduce the following parametrization of the set
$V(p,\e)$. This program has been proved, for the Laplacian
operator in \cite{BC}. The proof can be extended to our case
without any change. We consider the following minimization problem
for a function $u\in V(p,\e)$ with $\e$ small
\begin{eqnarray}\label{e:25}
\min \{\mid\mid
u-\sum_{i=1}^p\a_iP\d_{(a_i,\l_i)}\mid\mid_{2},\quad
\a_i>0,\quad\l_i>0,\quad a_i\in \O\}
\end{eqnarray}
\begin{pro}\label{p:22}
For any $p\in \N^*$, there exists $\e_p>0$ such that, for any
$0<\e<\e_p$, $u\in V(p,\e)$, the minimization problem \eqref{e:25} has a unique
solution (up to permutation). In particular, we can write 
$u\in V(p,\e )$ as follows 
$u=\sum_{i=1}^p\a_iP\d_{(a_i,\l_i)} + v$, where $(\bar{\a }_1,...,\bar{\a }_p,\bar{a}_1,...,\bar{a}_p,\bar{\l }_1,...,
\bar{\l }_p)$ is the solution of \eqref{e:25} and  $v$ satisfies
\begin{eqnarray*}
(V_0)\quad (v,P\d_{(a_i,\l_i)})_2=(v,\partial P\d_{(a_i,
\l_i)}/\partial \l_i)_2=0,\,  (v,\partial P\d_{(a_i,
\l_i)}/\partial a_i)_2=0 \mbox{ for }i=1,...,p,
\end{eqnarray*}
with $(u,w)_2=\int_{\O}\D u\D w$.
\end{pro}
Now we  recall the expansion of the functional
in the set $V(p,\e)$. To simplify notation, we will
write $\d_i$ instead of $\d_{(a_i,\l_i)}$.
\begin{pro}\label{p:23}\cite{BH}
There exists $\e_0>0$ such that for any
$u=\sum_{i=1}^p\alpha_iP\delta_i+v\in V(p,\e)$, $\e<\e_0$, $v$
satisfying $(V_0)$, we have
\begin{align*}
J(u)= & \frac{S_6^{2/3}\sum_{i=1}^p\a_i ^2}{(\sum_{i=1}^p\a_i ^6
K(a_i))^{1/3}}\biggl[1+\frac{c_0^6w_5}{24S_6\sum_{i =1}^p
K(a_i)^{-1/2}}\left(
\sum_{i=1}^p \frac{-\D K(a_i)}{60K(a_i)^{\frac{3}{2}}\l_i ^2}\right.\\
 & \left.+ \sum_{i=1}^p\frac{H(a_i,a_i)}{K(a_i)^{\frac{1}{2}}\l_i ^2}
  -\sum_{i\ne j}\frac{1}{(K(a_i)K(a_j))^{\frac{1}{4}}}\biggl(
\e_{ij}-\frac{H(a_i,a_j)}{\l_i\l_j} \biggr)\right) -f(v)\\
 & +\frac{1}{\sum_{i=1}^p\a_i ^2S_6}Q(v,v)+
o\left(\sum\frac{1}{\l_i ^2}+\frac{1}{(\l_id_i)^2}+\sum\e_{kr}+
\mid\mid v\mid\mid_{2}^2\right)\biggr],
\end{align*}
where
$$Q(v,v)=\mid\mid v\mid\mid_{ 2}^2-5
\frac{\sum_{i=1}^p\a_i ^2}{\sum_{i=1}^p\a_i ^6 K(a_i)}\int_\O K
\bigl(\sum_{i=1}^p\a_i P\d_i\bigr)^4v^2\quad ,
$$
$$f(v)=\frac{2}{\sum_{i=1}^p\a_i ^6 K(a_i)S_6}\int_\O
K \bigl(\sum_{i=1}^p\a_i P\d_i \bigr)^5v\quad ,
\qquad $$
$c_0$ is defined in \eqref{d}, $w_5$ is the volume of the unit
sphere $S^5$ and $S_6=\int_{\R^6}\frac{dx}{(1+|x|^2)^6}$.\\
Furthermore, if each $a_i$ is near a critical point $y_{j_i}$ of $K$ with $y_{j_i}\neq y_{j_k}$ for $i\neq k$, then this expansion
becomes
 \begin{align*}
 J(u) & = \frac{S_6^{2/3}\sum_{i=1}^p\a_i ^2}{(\sum_{i=1}^p\a_i ^6
K(a_i))^{\frac{1}{3}}}\biggl[1+\frac{c_0^6w_5}{24S_6
\sum_{i=1}^pK(a_i)^{-\frac{1}{2}}}\biggl(
\sum_{i=1}^p\left(\frac{-\D K(a_i)}{60K(a_i)^{3/2}}
+\frac{H(a_i,a_i)}{K(a_i)^{1/2}}\right)\frac{1}{\l_i ^2}\\
 & -\sum_{i\ne j}\frac{G(a_i,a_j)}{(K(a_i)K(a_j))^{1/4}}
\frac{1}{\l_i\l_j}\biggr)-f(v) + \frac{1}{\sum_{i=1}^p\a_i
^2S_6}Q(v,v)+o\biggl(\sum \frac{1}{\l_i ^2}+\mid\mid v\mid\mid_{2
}^2 \biggr)\biggr].
\end{align*}
\end{pro}
Now we are going to study the $v$-part of $u$. Let us observe that
\begin{eqnarray}\label{e:g}
\int K(\sum_{i=1}^p\a_iP\d_i)^4v^2=\sum_{i=1}^p
\a_i ^4K(a_i)\int P\d_i ^4v^2+o(\mid\mid v\mid\mid_{2}^2).
\end{eqnarray}
 Thus, using \eqref{e:g} and the fact that $\a_i ^4K(a_i)/(\a_j^4K(a_j))=1+o(1)$, we derive that the quadratic form $Q(v,v)$, defined in Proposition \ref{p:23},  is equal to
 $$\mid\mid v\mid\mid_{2}^2-5\sum_{i=1}^p\int
 P\d_i ^4v^2+o(\mid\mid v\mid\mid_{2}^2).$$
 Therefore, as in \cite{B1}, $Q(v,v)$ is a quadratic form positive
definite (see \cite{BE1}). It follows the following proposition

\begin{pro}\label{p:24}
There exists a $C^1$-map which, to each $(\a_1,...,\a_p, a_1,...,a_p, a, \l,...,\l_p)$ satisfying
$\sum_{i=1}^p\a_iP\d_{(a_i,\l_i)}\in V(p,\e)$, with $\e$ small
enough, associates $\ov{v}=\ov{v}(\a,a,\l)$ satisfying $(V_0)$
such that $\ov{v}$ is unique and minimizes \\
$J(\sum_{i=1}^p\a_iP\d_{(a_i,\l_i)}+v)$ with respect to $v$
satisfying $(V_0)$. Moreover, we have the following estimate
\begin{eqnarray*}
\mid\mid\ov{v}\mid\mid_{2} \leq c\mid f\mid
=O\biggl(\sum_{i=1}^p \frac{\mid\n
K(a_i)\mid}{\l_i}+\frac{1}{\l_i ^2}
+\sum{\e_{ij}(log\e_{ij}^{-1})^{1/3}}
+\frac{1}{(\l_id_i)^2}\biggr)
\end{eqnarray*}
\end{pro}
\section{Proof of Theorems}
\mbox{}
Before  giving the proof of our results, we first extract from \cite{BH} the characterization of   the critical points at infinity of our problem. We recall that the critical points at infinity are  the orbits of the gradient flow of $J$ which remain in $V(p,\e(s))$, where $\e(s)$ , a given function, tends to zero when $s$ tends to $+\infty$ (see \cite{B1}).
 \begin{pro}\label{p:31}  \cite{BH}
 Assume that for any  $s$, $M(\tau_s)$ is nondegenerate. Thus, there exists a pseudogradient $\ W\ $ so that the
following holds.\hfill\break There is a constant $c>0$ independent
of $ u=\sum_{i=1}^p\alpha_iP\delta_i$ in $V(p,\varepsilon)$ so
that
$$\bigl(-\partial J( u),W\bigr)\geq
c\sum_i\left(\frac{1}{\l_i ^2}+
\frac{1}{(\lambda_id_i)^3}+\sum_{k\ne
i}\varepsilon_{ki}^{3/2}+\frac{\mid\n K(a_i)\mid}{\l_i}
\right)\leqno{i)}$$
$$\biggl(-\partial J(u+\overline{v}),
W+\frac{\partial\overline{v}}{\partial(\alpha,a,\lambda)}(W)\biggr)
\geq c\sum_i\left(\frac{1}{\l_i ^2}+
\frac{1}{(\lambda_id_i)^3}+\sum_{k\ne
i}\varepsilon_{ki}^{3/2}+\frac{\mid\n K(a_i)\mid}{\l_i}
\right)\leqno{ii)} $$
 iii) $\displaystyle{\mid W\mid}$ is bounded.\hfill\break
 iv) $W$ satisfies the Palais-Smale condition
away from the critical points at infinity.\hfill\break
 v) The minimal distance to the boundary only
increases if it is small enough.\hfill\break
 vi) The $\lambda_i$'s are bounded away from
the case where each $a_i$ is near a critical point $y_{j_i}$
satisfying $j_i\ne j_k$ for $i\ne k$ and $\rho(y_{j_1},...,y_{j_p})>0$.
\end{pro}
\begin{cor}\label{c:32} \cite{BH}
Assume that  for any $s$, $M(\tau_s)$ is nondegenerate, and assume further that $J$ has no critical point in $\Sig^+$. Then the only critical
 points at infinity in $V(p,\e)$ correspond to $\sum_{j=1}^p
K(y_{i_j})^{-1/4}P\d_{(y_{i_j},\infty)}$, where $p\in N^*$ and the
points $y_{i_j}$'s are critical points of $K$ satisfying
$\rho(y_{i_1},...,y_{i_p}) > 0$.\\
 In addition, in the neighborhood of such a critical point at infinity, we can find a change of variable
$$
(a_1,...,a_p,\l_1,...,\l_p)\to ( a'_1,..., a'_p,\l '_1,...,\l '_p):=( a',\l ')
$$  
such that
\begin{eqnarray*}
J(\sum_{i=1}^p\alpha_iP\delta_{(a_i,\lambda_i)})=
\Psi(\a,a',\L'):= \frac{S_6^{2/3}\sum_{i=1}^p\alpha_i ^2
}{(\sum_{i=1}^p\alpha_i ^6K(a_i'))^{1/3}}\biggl(1+(c_1'-\eta)
{}^T\L' M(Y)\L' \biggr)
\end{eqnarray*}
where $\a=(\a_1,...,\a_p)$, $a'=(a_1',...,a_p')$ ,
${}^T\L'=(\l_1',...,\l_p')$, $c_1'$ is a positive constant and $\eta$ is small positive constant. 
\end{cor}
Now we are ready to prove our theorems.\\
\begin{pfn}{\bf Theorem \ref{t:11} }
For $\eta > 0$ small enough, we introduce, following \cite{BE1}, this neighborhood of $\Sig^+$ $$
V_\eta(\Sig^+)=\{u\in \Sig/\, \,e^{2J(u)} J(u)^2|u^-|_{L^6}^4 < \eta \},
$$
where $u^-=\max(0,-u)$.\\
Recall that, from Proposition \ref{p:31} we have a vector field $W$ defined in $V(p,\e)$ for $p\geq 1$. Outside $\cup_{p\geq 1}V(p,\e/2)$, we will use $-\n J$ and our global vector field $Z$ will be built using a convex combination of $W$ and $-\n J$. $V_\eta(\Sig^+)$ is invariant under the flow line generated by $Z$ (see \cite{BE1}).
 Arguing by contradiction, we assume that $J$ has no critical points in $V_\eta(\Sig^+)$. For any $y$ critical point of $K$, set
$$
c_\infty (y) = \biggl (\frac{S_6}{K(y)^{1/2}} \biggr )^{2/3}.
$$
Since $y_0$ is the unique absolute maximum of $K$, we derive that
$$
c_\infty (y_0) < c_\infty (y), \qquad \forall y \neq y_0,
$$
where $y$ is any critical point of $K$.\\
 Let $u_0 \in \Sig^+$ such that 
\begin{eqnarray}\label{e:31}
c_\infty (y_0) < J(u_0) < \inf _{y/ y\neq y_0, \n K(y) =0}c_\infty (y)
\end{eqnarray}
and let $\eta (s,u_0)$ be the one parameter group generated by $Z$. It is known that $|\n J|$ is lower bounded outside $V(p,\e /2)$, for any $p\in \N^*$ and for $\e$ small enough, by a fixed constant which depends only on $\e$. Thus the flow line $\eta (s,u_0)$ cannot remain outside of the set  $V(p,\e /2)$. Furthermore, if the flow line travels from  $V(p,\e /2)$ to the boundary of  $V(p,\e )$, $J(\eta (s,u_0))$ will decrease by a fixed constant which depends on $\e$. Then, this travel cannot be repeated in an infinite time. Thus there exist $p_0$ and $s_0$ such that the flow line enters into  $V(p_0,\e /2)$ and it does not exit from $V(p_0,\e)$. Since $u_0$ satisfies \eqref{e:31}, we derive that $p_0=1$, thus, for $s\geq s_0$,
$$
\eta (s,u_0)=\a _1 P\d_{ (x_1(s),\l _1(s))} + v(s)
$$
Using again \eqref{e:31}, we deduce that $x_1(s)$ is outside $\mathcal{V}(y,\tau )$ for any $y \in \mathcal{F}^+ \diagdown \{y_0\}$, where $\mathcal{V}(y,\tau )$ is a neighborhood of $y$ and where $\tau$ is a small positive real. Now by  assumptions of Theorem \ref{t:11} and by the construction of a pseudogradient $Z$, we derive that $\l _1(s)$ remains bouded along the flow lines of $Z$. Thus, we obtain
$$
|\n J(\eta (s,u_0)).Z(\eta (s,u_0))| \geq c > 0 \quad \forall s\ge 0,
$$
where $c$ depends only on $u_0$.\\
 Then when $s$ goes to $+\infty$, $J(\eta (s,u_0))$ goes to $-\infty$ and this yields a contradiction.
Thus there exists a critical point of $J$ in $V_\eta(\Sig^+)$. Arguing as in  Proposition 4.1 \cite{BH}(see also Proof of Theorem 1.1 of \cite{BE1}), we prove that such a critical point is positive and hence our result follows.
\end{pfn}\\
\begin{pfn}{\bf Theorem \ref{t:12} }
Arguing by contradiction, we assume that $(1)$ has no solution. Notice that under the assumption of our theorem, we have $\rho (y_i,y_j) < 0$ for $y_i\neq y_j \in \mathcal{F}^+$. Thus using Corollary \ref{c:32} we derive that the only critical points at infinity of $J$ are in one to one correspondance with the critical points of $K$ such that $y_i\in\mathcal{F}^+$. The unstable manifold $W_u(y_i)_\infty$ of such a critical point at infinity can be described using Corollary \ref{c:32} as the product of $W_s(y_i)$ (for a pseudo gradient of $K$) by $[A,+\infty ]$ domaine of the variable $\l$, for some postive real $A$ large enough. \
Since $J$ has no critical points in $\Sig^+$, it follows that $\Sig^+$ retracts by deformation on $\bigcup_{y_i\in\mathcal{F}}W_u(y_i)_\infty :=X_\infty$ (see Sections $7$ and $8$ of \cite{BR}). Now, using the fact that $\Sig^+$ is a contractible set, we derive that $X_\infty$ is contractible leading to the contractibility of $X$, which is in contradiction with our assumption. Hence our result follows.
\end{pfn}\\
Now before giving the proof of Theorem \ref{t:13}, we state the following lemma. Its proof is very similar to the proof of Corollary B.3 of \cite{BC}(see also \cite{B2}), so we will omit it.
\begin{lem}\label{l:33}
Let $a_1$, $a_2 \in \O$, $\a_1$, $\a_2 > 0$ and $\l$ large enough. For  $u=\a_1 P\d_{ (a_1,\l )} +\a_2P\d_{ (a_2,\l )}$, we have 
$$
J\biggl(\frac{u}{||u||}\biggr) \leq S_6^{2/3}\biggl(\frac{1}{K(a_1)^{1/2}} + \frac{1}{K(a_2)^{1/2}} \biggr)^{1/3}(1+o(1)).
$$
\end{lem}
\begin{pfn}{\bf Theorem \ref{t:13} }
Again, we argue by contradiction. We assume that $J$ has no critical point in $V_\eta (\Sig^+)$. Let
$$
c_\infty (y_0,y_{i_1}) = S_6^{2/3}\biggl(\frac{1}{K(y_0)^{1/2}} + \frac{1}{K(y_{i_1})^{1/2}} \biggr)^{1/3}
$$
 We observe that under the assumption $(A_1)$ of Theorem \ref{t:13}, $(y_0,y_{i_1})$ is not a critical point at infinity of $J$. Using Corollary \ref{c:32} and the assumption $(A_2)$ of Theorem \ref{t:13}, it follows that the only critical points at infinity of $J$ under the level $c_1=c_\infty (y_0,y_{i_1}) + \e$, for $\e$ small enough, are $P\d_{ (y_0,\infty )}$ and   $P\d_{ (y_{i_1},\infty )}$. 
 The unstable manifolds at infinity of such  critical
points at infinity, $W_u(y_0)_\infty$, $W_u(y_{i_l})_\infty$ can be
described, using Corollary \ref{c:32}, as the product of $W_s(y_0)$,
 $W_s(y_{i_l})$ (for a pseudogradient of $K$ ) by $[A, +\infty [$
 domaine of the variable $\l$, for some positive
number $A$ large enough.\\
Since $J$ has no critical point, it follows that $ J_{c_1}=\{u\in
\sum ^+ / J(u) \leq c_1 \}$ retracts by deformation on $X_\infty =
 W_u(y_0)_\infty \cup W_u(y_{i_1})_\infty$ (see Sections 7 and 8 of
\cite{BR}) which can be parametrized  by $X \times
[A, +\infty[$, where $X=\overline{W_s(y_{i_1})}$.\\
Under ${\bf 1.3}$ and ${\bf 1.4}$ (see first section), we have $X= W_s(y_0) \cup W_s(y_{i_1})$. Thus $X$ is a manifold in dimension $k_1$ without boundary.\\
Now we claim that $X_\infty$ is contractible in $J_{c_1}$. Indeed, let $h:  [0,1] \times X \times \left[0,\right.+\infty\left[ \right.\longmapsto  \Sigma ^+$ defined by
$$
 (t,x,\l _1 ) \longmapsto \frac{tP\d _{(y_0,\l )} + (1-t)P\d _{(x,\l )}}{||tP\d _{(y_0,\l )} + (1-t)P\d _{(x,\l )}||}
$$
$h$ is continuous and satisfies 
$$
h(0,x,\l )=\frac{P\d _{(x,\l )}}{||P\d _{(x,\l )}||} \quad \mbox{and}\quad  h(1,x,\l )=\frac{P\d _{(y_0,\l )}}{||P\d _{(y_0,\l )}||}.
$$
In addition, since $K(x)\geq K(y_{i_1})$ for any $x\in X$, it follows from Lemma \ref{l:33} that $J(h(t,x,\l ))< c_1$, for any $(t,x,\l )\in [0,1]\times X\times[A,+\infty )$. Thus the contraction $h$ is performed under the level $c_1$.We derive that 
  $X_\infty$ is contractible in $J_{c_1}$, which retracts by deformation on $X_\infty$,
therfore $X_\infty$ is contractible leading to the contractibility of
$X$, which is a contradiction since $X$ is a manifold in dimension $k_1$ without boundary. Hence  there exists a critical point of $J$ in $V_\eta (\Sigma ^+ )$. Arguing as in Proposition 4.1 \cite{BH}, we prove that such a critical point is positive. Now, we are going to show that such a critical point has a Morse index equal to $k_1$ or $k_1+1$.\\
Using a dimension argument and since $h([0,1], X_\infty )$ is a manifold in dimension $k_1+1$, we derive that the Morse index of such a critical point is $\leq k_1+1$.\\
Now, arguing by contradiction, we assume that the Morse index is $\leq k_1-1$. 
Perturbing, if necessary $J$, we may assume that all the critical
points of $J$ are nondegenerate and have their Morse index $\leq
k_1-1$.\\
 Such critical points do not change the homological group in
dimension $k_1$ of level sets of $J$.\\
Now, let  $c_\infty (y_{i_1})= S_6^{2/3}K(y_{i_1})^{-1/3}$ and let $\e$ be a small positive real.
Since $X_\infty$ defines a homological class in dimension $k_1$ which is
trivial in  $J_{c_1}$, but nontrivial in  $J_{c_\infty (y_{i_1})+ \e}$, our result follows.
\end{pfn}\\
\begin{pfn}{\bf Theorem \ref{t:14}}
We notice that the assumption $(H_4)$ implies that $(y_0,y_{i_1})$ is a critical point at infinity of $J$. Now, arguing by contradiction, we assume that $(1)$ has no solution. We claim that $f_\l (C_{y_o}(X))$ retracts by deformation on $X\cup W_u(y_0,y_{i_1})_\infty$. Indeed, let 
$$
u=\a P\d _{(y_0,\l )} + (1-\a )P\d _{(x,\l )} \in f_\l (C_{y_0}(X)),
$$
the action of the flow of the pseudo gradient $Z$ defined in the proof of Theorem \ref{t:11} is essentially on $\a$.\\
- If $\a < 1/2$, the flow of $Z$ brings $\a$ to zero and thus $u$ goes in this case to $\overline{W_u(y_0)}_\infty \equiv \{y_0\}$.\\
- If  $\a > 1/2$, the flow of $Z$ brings $\a$ to $1$ and thus $u$ goes, in this case, to $\overline{W_u(y_{i_1})}_\infty \equiv X_\infty$.\\
- If $\a = 1/2$, since only $x$ can move then $y_0$ remains one of the points of concentration of $u$ and $x$ goes to $W_s(y_i)$, where $y_i=y_{i_1}$ or $y_i=y_0$ and two cases may occur :\\
- In the first case, that is, $y_i=y_{i_1}$, $u$ goes to $W_u(y_0, y_{i_1})_\infty$ .\\
- In the second case, that is,$y_i=y_0$, there exists $s_0 \geq 0$ such that $x(s_0)$ is close to $y_0$.  Thus, using Lemma \ref{l:33}, we have the following inequality
$$
J(u(s_0)) \leq c_\infty (y_0,y_0) + \gamma := c_2,
$$
where $c_\infty (y_0,y_0)= S_6^{2/3} (\frac{2}{K(y_0)^{1/2}})^{1/3}$ and where $\gamma$ is a positive constant small enough.\\
Now, using  assumption $(H_5)$, it follows from Corollary \ref{c:32} that $J_{c_2}$ retracts by deformation on $W_u(y_0)_\infty \equiv \{y_0\}$ and thus $u$ goes to $W_u(y_0)_\infty$. Therefore $f_\l (C_{y_0}(X))$ retracts by deformation on $X_\infty \cup W_u(y_0,y_{i_1})_\infty$.\\
Now, since $\mu (y_{i_1}) =0$, it follows that this strong
retract  does not intersect
$W_u(y_0,y_{i_1})_\infty$ and thus it is contained in $X_\infty$. Therefore $X_\infty$ is contractible, leading to the contractibility of $X$, which is a contradiction since $X$ is a manifold of dimension $k_1$ without boundary.  Hence $(1)$ admits a solution. Now, using the same arguments as those used in the proof of Theorem \ref{t:13}, we easily derive that the Morse index of the solution provided above is equal to $k_1$ or $k_1+1$. Thus our result follows.
\end{pfn}\\
\begin{pfn}{\bf Theorem \ref{t:15}}
Assume that $(1)$ has no solution. Using the same arguments as those in the proof of Theorem \ref{t:14}, we deduce that  $f_\l (C_{y_0}(X))$ retracts by deformation on $$
X_\infty \cup  \left(\cup
_{y\in B_k}W_u(y_0, y)_\infty \right)\cup D,
$$
 where $D\subset \sigma$ is a stratified set and where
$\sigma = \cup _{y\in X\diagdown B_k}W_u(y_0, y)_\infty$
is a manifold in dimension at most $k$.\\
Since $\mu  (y)=0$ for each $y \in B_k$,  $f_\l (C_{y_0}(X))$
retracts by deformation on $X_\infty \cup D$, and therefore $H_*(X_\infty\cup D)=0$,
for all $*\in \N^*$, since $f_\l(C_{y_0}(X))$ is a contractible set. Using the exact homology sequence of $(X_\infty\cup D,X_\infty)$, we obtain 
$$
...\to H_{k+1}(X_\infty\cup D) \to ^{\pi}  H_{k+1}(X_\infty\cup D, X_\infty) \to
^{\partial}  H_k(X_\infty) \to ^{i}  H_k(X_\infty\cup D) \to ...
$$
Since $H_*(X_\infty\cup D) = 0$, for all $*\in \N^*$, then
$H_k(X_\infty)=H_{k+1}(X_\infty\cup D, X_\infty)$.\\
In addition, $(X_\infty\cup D, X_\infty)$ is a stratified set of dimension at most
$k$, then $H_{k+1}(X_\infty\cup D,X_\infty) = 0$, and therefore $H_k(X_\infty)=0$.  This implies that $H_k(X)=0$
( recall that $X_\infty\equiv X\times [A,\infty)$). This
yields a contradiction since $X$ is a manifold in dimension $k$
without boundary. Then, arguing as in the end of the proof of Theorem \ref{t:14}, our theorem follows.
\end{pfn}\\


\begin{thebibliography}{99}

\bibitem{B1}
 A. Bahri,
  \emph{ Critical points at infinity in some variational problems,}
 Pitman Res.  Notes Math. Ser. \textbf{182}, Longman Sci. Tech. Harlow (1989).
\bibitem{B2}
 A. Bahri,
  \emph{ An invariant for Yamabe-type flows with application to
  scalar curvature problems in high dimension, } A celebration of
  J.F. Nach Jr., Duke Math. J. \textbf{81} (1996), 323-466.
\bibitem{BC}
 A. Bahri and J.M. Coron,
\emph{ On a nonlinear Elliptic equation Involving the critical
Sobolev Exponent: The effect of the topology on the Domain,}
 Comm. Pure Appl. Math. \textbf{41} (1988), 253-294.
\bibitem{BaC}
 A. Bahri and J. M. Coron,
\emph{The scalar curvature problem on the standard three dimensional spheres,}
 J. Funct. Anal. \textbf{95} (1991), 106-172.
\bibitem{BR}
 A. Bahri and P.H. Rabinowitz,
 \emph{ Periodic solution of Hamiltonian Systems of three-body type,
 }Ann. Inst. Henri Poincar\'e Anal. Non lin\'eaire, \textbf{8}
 (1991), 561-649.
\bibitem{BCCH}
 M. Ben Ayed, Y. Chen, H. Chtioui and M. Hammami,
\emph{ On the prescribed scalar curvature problem on 4-manifolds,}
 Duke Math. J. \textbf{84} (1996), 633-677.
\bibitem{BCH}
 M. Ben Ayed, H. Chtioui and M. Hammami,
\emph{ The scalar curvature problem on higher dimensional spheres,}
 Duke Math. J. \textbf{93} (1998), 379-424.
\bibitem{BE1}
 M. Ben Ayed and K. El Mehdi,
 \emph{ The Paneitz Curvature problem on lower dimensional spheres,} Preprint The Abdus Salam ICTP, Trieste, Italy, IC/2003/48.
\bibitem{BE2}
 M. Ben Ayed and K. El Mehdi,
 \emph{ Existence of conformal metrics on spheres with prescribed  Paneitz Curvature,} Preprint ICTP The Abdus Salam ICTP, Trieste, Italy, IC/2003/55.
\bibitem{BEH}
M. Ben Ayed, K. El Mehdi and M. Hammami,
\emph{ Some existence results for a fourth order elliptic equation involving critical exponent, } Preprint  The Abdus Salam ICTP, Trieste, Italy, IC/2003/49.
\bibitem{BH}
M. Ben Ayed and M. Hammami,
\emph{On a fourth order elliptic equation with critical nonlinearity in dimension six, } Preprint 2003.
\bibitem{BGP}
F. Bernis, J. Garcia-Azorero and I. Peral,
\emph{ Existence and multiplicity of nontrivial solutions in semilinear critical problems,} Adv. Differential Equations \textbf{1} (1996), 219-240.
\bibitem{BrC}
H. Brezis and J.M. Coron,
\emph{Convergence of solutions of $H$-systems or how to blow bubbles,}
Arch. Rational Mech. Anal. \textbf{89} (1985), 21-56.
\bibitem{C}
H. Chtioui,
\emph{Presribing the scalar curvature problem on three and four manifolds,} 
Preprint \textbf{} (2003).
\bibitem{CE}
H. Chtioui and K. El Mehdi,
\emph{Prescribed scalar curvature with minimal boundary mean curvature on $S^4_+$, } Preprint  The Abdus Salam ICTP, Trieste, Italy, IC/2003/57.
\bibitem{DHL}
 Z. Djadli, E. Hebey and M. Ledoux,
 \emph{ Paneitz type operators and applications, } Duke Math. J. \textbf{104} (2000), 129-169.
\bibitem{DMO}
 Z. Djadli, A. Malchioldi and M.Ould Ahmedou,
 \emph{ Prescribing a fourth order conformal invariant on the
standard sphere- Part I: a perturbation result, }Comm. Contemp. Math. \textbf{4} (2002), 375-408, \emph{ Part II: blow-up analysis and applications, } Annali Scuola Norm. Sup. Pisa \textbf{5} (2002), 387-334.
\bibitem{EFJ}
D.E. Edmunds, D. Fortunato and E. Janelli,
\emph{ Critical exponents, critical dimension and the biharmonic operator, }
Arch. Rational Mech. Anal. \textbf{112} (1990), 269-289.
\bibitem{EO}
 F. Ebobisse and M. Ould Ahmedou,
 \emph{ On a nonlinear fourth order elliptic equation involving
the critical Sobolev exponent, } Nonlinear Anal. TMA \textbf{52} (2003), 1535-1552.
\bibitem{F}
V. Felli,
\emph{Existence of conformal metrics on $S^n$ with prescribed
  fourth-order invariant,}
Adv. Differential Equations \textbf{7} (2002), 47-76.
\bibitem{HR}
 E. Hebey and F. Robert,
 \emph{ Coercivity and Struwe's compactness for Paneitz type
operators with constant coefficients, } Calc. Var. Partial Differential Equations \textbf{13} (2001), 491-517.
\bibitem{HV}
J. Hulshof and R.C.A.M. Van Der Vorst,
\emph{ Differential systems with strongly indefinite variational structure, } J. Funct. Anal. \textbf{114} (1993), 32-58.
\bibitem{Lin}
C.S. Lin,
\emph{Classification of solutions of a conformally invariant fourth order equation in $\R^n$,} 
Comment. Math. Helv. \textbf{73} (1998), 206-231.
\bibitem{Lio}
P. L. Lions,
\emph{ The concentration compactness principle in the calculus of variation. The limit case,} Rev. Mat. Iberoamaricana \textbf{1} (1985), I: 165-201; II: 45-121.
\bibitem{M}
J. Milnor,
\emph{Lecturess on h-Cobordism Theorem,}
 Princeton University Press, Princeton \textbf{} 1965.
\bibitem{NSJ}
E.S. Noussair, C.A. Swanson and Y. Jianfu,
\emph{Critical semilinear biharmonic equations in $\R^n$, } Proc. Royal Soc. Eidinburg \textbf{121A} (1992), 139-148.
\bibitem{PV}
L. Peletier and R.C.A.M. Van Der Vorst,
\emph{ Existence and nonexistence of positive solutions of nonlinear elliptic systems and the biharmonic equation, } Differential Integral Equations \textbf{5} (1992), 747-767. 
\bibitem{PS}
P. Pucci and J. Serrin,
\emph{Critical exponents and critical dimensions for polyharmonic operator,} J. Math. Pures Appl. \textbf{69} (1990), 55-83.
\bibitem{R}
O. Rey,
\emph{Bifurcation from infinity in a nonlinear elliptic equation involving the limiting Sobolev exponent,} Duke Math. J. \textbf{60} (1990), 815--861. 
\bibitem{S}
  M. Struwe,
\emph{ A global compactness result for elliptic boundary value
problem involving limiting nonlinearities}, Math. Z. \textbf{187}
 (1984), 511-517.
\bibitem{V1}
  R.C.A.M. Van Der Vorst,
\emph{ Fourth order elliptic equations with critical growth, } C.
R. Acad. Sci. Paris, t. 320, S\'erie I, (1995), 295-299.
\bibitem{V2}
  R.C.A.M. Van Der Vorst,
\emph{ Best constant for the embedding of the space $H^2\cap H_0^1
(\O)$ into $L^{2n/(n-4)}(\O)$, } Diff. Int. Equa. \textbf{6}
(1993), 259-276.

\end{thebibliography}
\end{document}